\documentclass[12pt]{amsart}

\newtheorem{theorem}{Theorem}[section]

\newtheorem{remark}[theorem]{Remark}
\newtheorem{corollary}[theorem]{Corollary}

\begin{document}
\title[Diagonals in Tensor Products of Operator Algebras]
{Diagonals in Tensor Products of Operator Algebras}

\vspace{30 mm}

\author{Vern I. Paulsen and Roger R. Smith}
\address{Department of Mathematics\\University of Houston\\4800 Calhoun
Road\\
Houston, TX 77204--3476 U.S.A.}

\email{vern@math.uh.edu}
\address{Department of Mathematics\\Texas A$\&$M University\\
College Station, TX 77843--3368 U.S.A. }
\email{rsmith@math.tamu.edu}
\date{September 30, 2000.}
\thanks{Both authors were partially supported by grants from the NSF. The
first author also wishes to thank
the Department of Mathematics, Rice University where parts of this research
were completed.}

\maketitle

\vspace{40 mm}

\begin{abstract}
In this paper we give a short, direct proof, using only properties of the
Haagerup tensor
product, that if an operator algebra $A$ possesses
a diagonal in the Haagerup tensor product of $A$ with itself, then $A$ must
be isomorphic
to a finite dimensional $C^*$-algebra. Consequently, for operator algebras, the
first
Hochschild cohomology group,
$H^1(A,X) = 0$ for every bounded, Banach $A$-bimodule $X$, if and only if
$A$ is isomorphic
to a finite dimensional $C^*$-algebra.
\end{abstract}


\pagebreak
\newpage
\baselineskip = 29pt
\section{Introduction.}
Let $A$ be a complex algebra with unit 1. A {\em diagonal in $A \otimes A$}
is an element
$u = \sum a_i \otimes b_i$ such that $\sum a_i b_i = 1$ and
$\sum (aa_i) \otimes b_i = \sum a_i \otimes (b_ia)$ for every $a \in A$.
For example, if $\mathbb{M}_n$ denotes the algebra of $n \times n$ complex matrices and
$E_{ij}$
denotes the standard matrix units, then
$\sum_{i=1}^{n} E_{i1} \otimes E_{1i}$ is easily seen to be a diagonal in $\mathbb{M}_n
\otimes \mathbb{M}_n$.

It is fairly well-known that the existence of a diagonal is equivalent to
the vanishing of the first Hochschild cohomology, $H^1(A,X)$
for every $A$-bimodule $X$.
Since this fact is elementary, we quickly recall a proof, for clarity.
First, assume that we have a diagonal $u$, as above, and that we are given
an $A$-bimodule $X$
and a derivation $\delta : A \rightarrow X.$
If we set $x = \sum \delta (a_i)b_i,$ then it is easily checked that $\delta
(a) = xa - ax$ and
so every derivation into $X$ is inner.  That is, $H^1(A,X) = 0$, for every
$X.$
To prove the converse, one simply considers the $A$-bimodule $A \otimes A$
and
lets $X$ be the submodule which is the kernel of the product map. The map
$\delta : A \rightarrow X$
given by $\delta (a) = a \otimes 1 - 1 \otimes a$ is easily seen to be a
derivation.
If $w \in X$ is the element that implements this necessarily inner
derivation, then is is easily
checked that $u = 1 \otimes 1 - w$ is the desired diagonal.

The above proof easily extends to the case of various topological algebras,
where the module actions
and derivations are restricted to those which are continuous in some
appropriate sense. The only
change that must be made is that the algbraic tensor product of $A$ with
itself is replaced by its
completion in some appropriate topology.

In this setting, Helemskii, \cite{He}, and Selivanov, \cite{Se}, proved that a
$C^*$-algebra $A$
has the property that every bounded derivation into every bounded
$A$-bimodule is inner if and only
if $A$ is finite dimensional. By the above remarks this can be seen to be
equivalent to
characterizing those $C^*$-algebras $A$ which possess a diagonal in the
projective tensor product
of $A$ with itself.
The first author extended this result, \cite{Pa}, by proving that
a $C^*$-algebra $A$
has the property that every completely bounded derivation into every operator $A$-bimodule is inner if and only if $A$ is finite dimensional.
Again by the above remarks, this latter result is equivalent to proving that the only $C^*$-algebras $A$ which 
possess a diagonal in the Haagerup tensor product of $A$ with itself are the
finite dimensional
$C^*$-algebras.
Because the projective tensor norm is larger than the Haagerup tensor norm, this latter result implies the result of Helemskii and Selivanov.

Unfortunately, all of the proofs cited above relied on non-trivial results.
In particular, the proof in \cite{Pa} relied on deep results about nuclear and injective C*-algebras.
The purpose of this paper is to state a more general result, and to give a
short, self-contained proof using only properties of the
Haagerup tensor product.
Specifically, we show  that if $A$ is any algebra of operators on a Hilbert
space with a diagonal in the Haagerup
tensor product of $A$ with itself, then $A$ is necessarily isomorphic to a
finite direct sum of
matrix algebras.\newpage

\section{Main Results}

In this section we present our main results.
Let $H$ be a Hilbert space and let $B(H)$ denote the algebra of bounded
linear operators on $H.$
We let $A$ be any subalgebra of $B(H)$ which contains the identity operator.
We shall call each such algebra an {\em algebra of operators}, and we note
that we do not require it to be self-adjoint.

We briefly recall the definition of the Haagerup tensor product. Given $w
\in A \otimes A$
we set
\begin{equation}\label{eq2.1}\|w\|_h = {\mathrm {inf}}\ \left\{{ \left\| \sum
a_i a^{*}_i \right\|^{1/2}\,\left\| \sum b^{*}_i b_i \right\|^{1/2}
}\right\}\end{equation} where the infimum is taken over all ways to express $w$
as a finite sum  $\sum a_i \otimes b_i $ of elementary tensors.
This quantity defines a norm on $A \otimes A$ called the {\em Haagerup
tensor norm}
and the completion of $A \otimes A$ in this norm is called the {\em Haagerup
tensor product}
of $A$ with itself and is denoted $A \otimes_{h} A.$

This tensor norm has two very nice properties that we shall use.
The first is that any $w$ in the completion has a representation as a norm
convergent series,
$w = \sum^{\infty}_{i=1} a_i \otimes b_i$ with $\| \sum^{\infty}_{i=1} a_ia^{*}_i \|$
and $\| \sum^{\infty}_{i=1} b^{*}_ib_i \|$ both finite.
The second is that such a representation may be chosen so that $\{a_i\}^{\infty}_{i=1}$ and
$\{b_i\}^{\infty}_{i=1}$ are both
{\em strongly independent} sets in the following sense.
A sequence of elements $\{a_i\}^{\infty}_{i=1}$ which defines a bounded
operator $(a_1,a_2,\hdots)\in B(H^{\infty},H)$ is strongly independent if
the equation $\sum^{\infty}_{i=1}\lambda_i a_i=0$, where
$\{\lambda_i\}^{\infty}_{i=1} \in \ell^2$, can only be satisfied by
$\lambda_i=0$, $i \geq 1$. An equivalent formulation, \cite[Lemma 2.2]{ASS}, is
that the subspace \[\{(\phi(a_1),\phi(a_2),\hdots):\ \phi \in A^*\}\] is norm dense in
$\ell^2$. For these facts about the Haagerup tensor product we refer the
reader to \cite{BS,SS,Sm}.

\begin{theorem}\label{thm2.1} Let $A$ be an algebra of operators on a Hilbert
space.  If there is a diagonal in $A \otimes_{h} A$, then $A$ is finite
dimensional.
\end{theorem}

\noindent \begin{proof}
Let $u = \sum^{\infty}_{i=1} a_i \otimes b_i$ be a diagonal, where the series
is norm convergent and
$\{a_i\}^{\infty}_{i=1}$ and $\{b_i\}^{\infty}_{i=1}$ are strongly independent.
Since 
\begin{equation}\label{eq2.2}\sum^{\infty}_{i=1} a_ib_i = 1\end{equation}
is a norm convergent series, we may choose $M$ so that
\begin{equation}\label{eq2.3}\left\| \sum^{M}_{i=1} a_ib_i - 1 \right\| <
1/2,\end{equation}
 and we set $c = (\sum^{M}_{i=1} a_ib_i)^{-1}$. From the
Neumann series we know that $\|c\| < 2$. Now define two constants $k$ and
$\varepsilon$ by  
\begin{equation}\label{eq2.4} k=\mathrm{max}\ \left\{\left\|
\sum^{\infty}_{i=1} a_ia^{*}_i \right\|^{1/2},\ \ \left\| \sum^{\infty}_{i=1}
b^{*}_ib_i \right\|^{1/2}\right\},\ \ \ \varepsilon =
(8Mk^2)^{-1}.\end{equation}  
Since, for each $x \in A$, the series  
\begin{equation}\label{eq2.5}\sum^{\infty}_{i=1} xa_i \otimes b_i =
\sum^{\infty}_{i=1} a_i \otimes b_ix \end{equation}
are norm convergent, we may
apply, by \cite[Prop. 3.7]{BS}, an element $\phi \in A^*$ to (\ref{eq2.5}) to
obtain \begin{equation}\label{eq2.6}\sum^{\infty}_{i=1} \phi(xa_i)b_i =
\sum^{\infty}_{i=1} \phi(a_i)b_ix.\end{equation}
 From the strong independence of
$\{a_i\}^{\infty}_{i=1}$, we may choose linear functionals $\phi_j \in A^*$,
\mbox{$1 \leq j \leq M$,} such that 
\begin{equation}\label{2.7}\|(\phi_j(a_1), \phi_j(a_2),\hdots) - e_j
\|_2 < \varepsilon,\ \ \ 1 \leq j \leq M,\end{equation}
where $\{e_j\}^{\infty}_{j=1}$ denotes the
canonical orthonormal basis for $\ell^2$.

It now follows that 
\begin{equation}\label{eq2.8}\left\| b_jx - \sum^{\infty}_{i=1} \phi_j(a_i)b_ix
\right\| \leq \varepsilon \left\| \sum^{\infty}_{i=1} x^*b^*_ib_ix
\right\|^{1/2} \leq \varepsilon k \|x\|,\end{equation}
 for $1 \leq j \leq M,$ and for all $x
\in A$. Using (\ref{eq2.6}), we have that
\begin{equation}\label{eq2.9}\left\| b_jx - \sum^{\infty}_{i=1} \phi_j(xa_i)b_i
\right\| \leq \varepsilon k \|x\|,\end{equation}
 for $1 \leq j \leq M$,
and for all $x \in A$.

Since $\mathrm{lim}_{n\rightarrow\infty}\left\| \sum^{\infty}_{i=n} a_ia^*_i
\right\|=0$, we may
choose $N$ sufficiently large that
\begin{equation}\label{eq2.10}\left\| b_jx - \sum^{N}_{i=1} \phi_j(xa_i)b_i
\right\| \leq 2\varepsilon k \|x\|\end{equation}
holds for $1 \leq j \leq M$,
and for all $x \in A$.
The inequality $\|a_j\| \leq k$ follows from (\ref{eq2.4}), and so the
relation 
\begin{equation}\label{eq2.11}\left\| \sum^{M}_{j=1}\ \left[a_jb_jx -
\sum^{N}_{i=1} \phi_j(xa_i)a_jb_i\right] \right\| \leq 2\varepsilon Mk^2
\|x\|\end{equation}
is a consequence of multiplying the expression in (\ref{eq2.10}) on the left by
$a_j$ and summing over $j$. Now multiply (\ref{eq2.11}) on the left by $c$ and
use (\ref{eq2.4}) to obtain  
\begin{equation}\label{eq2.12}\left\|
x - \sum^{M}_{j=1} \sum^{N}_{i=1} \phi_j(xa_i)ca_jb_i \right\| \leq  4\varepsilon
Mk^2 \|x\| \leq \|x\|/2,\ \ \ x \in A.\end{equation}

Define a finite dimensional subspace of $A$ by
 \[B = \mathrm{Span}\,\{ ca_jb_i : 1
\leq j \leq M,\ 1 \leq i \leq N \}.\]
 The inequality (\ref{eq2.12}) implies that the Banach space
quotient map from $A$ to $A/B$ has norm at most $1/2$,
which can only happen when $A = B$.
We conclude that $A$ is finite dimensional.
\end{proof}

Ruan, \cite{Ru} has introduced another Hochschild cohomology for operator algebras which uses a family of maps called the {\em jointly completely bounded} maps. The relevant tensor norm for this cohomology is called the {\em operator space projective tensor norm.}

\begin{corollary} Let $A$ be an algebra of operators on a Hilbert space. 
If there is a diagonal in either the projective or operator space
projective tensor product of $A$ with itself, then $A$ is finite dimensional.
 \end{corollary}

 \noindent \begin{proof} For any element in the algebraic tensor product $A
\otimes A$, we have that its projective
 tensor norm is at least as large as its Haagerup norm. Thus, the identity map
on $A \otimes A$
 extends to a contractive map from the projective tensor product to the
Haagerup tensor product.
It is easily checked that if $u$ is a diagonal in the projective tensor
product, then its image under this map
 is a diagonal in the Haagerup tensor product, and the result follows from
Theorem \ref{thm2.1}. A similar argument applies to the operator space
projective tensor product.  \end{proof}

 Note that, since any $A$ as above is finite dimensional, the algebraic
tensor product is
 complete in every tensor norm. Thus we are reduced to the purely algebraic
 problem of determining those finite dimensional complex algebras $A$ that have
a diagonal in
 $A \otimes A$, which
  is essentially Burnside's theorem. We supply a simple proof below
that is based on the
 ideas that we have already introduced.

 \begin{theorem} Let $A$ be a finite dimensional, unital, complex algebra.
 If $A$ has a diagonal in $A \otimes A$, then $A$ is isomorphic to a direct
 sum of matrix algebras.
 \end{theorem}

 \noindent \begin{proof} Since $A$ can be represented as the algebra of left
multiplication operators on itself,
 we may assume that $A$ is a subalgebra of $\mathbb{M}_n$ for some $n$.

 Now suppose that $p$ is an invariant orthogonal projection for $A$, that
is, $pap = ap$ for all
 $a$ in $A$.
 Then it is easily seen that $X \equiv p\mathbb{M}_n(1-p)$ is an $A$-bimodule,
and that the equation
 \begin{equation}\label{eq2.13}\delta (a) \equiv pa(1-p) = pa -
ap\end{equation} defines a derivation of  $A$ into $X$.
 By hypothesis, there exists $x \in X$ such that 
\begin{equation}\label{eq2.14}\delta (a) = ax - xa,\ \ \ a \in A.\end{equation}
 Combining these last two equations, we see that $(p+x)$ commutes with $A$. 
Since $x^2 = 0$, the element $y \equiv 1 + x$ is
 invertible in $\mathbb{M}_n$ with inverse $y^{-1} = 1 - x$.

 The equations 
\begin{equation}\label{eq2.15} p(1+x)a(1-x) = (p+x)a(1-x)
= a(p+x)(1-x) = ap\end{equation}
 and 
\begin{equation}\label{eq2.16} (1+x)a(1-x)p = (1+x)ap = (1+x)pap = pap =
ap\end{equation}  show that $p$ commutes with $yAy^{-1}$, and thus
reduces this algebra.

 By inductively choosing such projections $p$ and conjugating by the
corresponding invertible elements,
 we may assume that the representation  $\pi : A \rightarrow \mathbb{M}_n$ is a
finite direct
 sum of representations, $\pi_i : A \rightarrow \mathbb{M}_{n_i},\ 
i=1,\hdots,k$, where the
 image $\pi_i(A)$ is a subalgebra of $\mathbb{M}_{n_i}$ that has no non-trivial
invariant projections.
Thus, $\pi_i(A)$ is a transitive subalgebra, and hence $\pi_i(A) =
\mathbb{M}_{n_i}$ by Burnside's theorem,
 \mbox{\cite[Cor. 8.6]{RR}.}

 Using the simplicity of each matrix algebra, it is now easy to argue that
$A$ is isomorphic
 to a finite direct sum of matrix algebras.
To see this, note that if $J_i = ker( \pi_i)$, then $\pi_j(J_i)$ is either
$\mathbb{M}_{n_j}$ or $(0)$ and
argue by induction on $k$.

\end{proof}

\begin{corollary} If $A$ is an algebra of operators and $A$ has
a diagonal in one of the
Haagerup, projective, or operator space projective tensor products of $A$ with
itself, then $A$ is isomorphic to a
finite direct sum of matrix algebras.
\end{corollary}

We end by formally stating the equivalent theorems in terms of Hochschild
cohomology.
If $A$ is any Banach algebra, then by an $A$-{\em bimodule} $X$ we mean
any Banach
space $X$ equipped with an $A$-bimodule action satisfying, $\|axb\| \leq
c\|a\| \|x\| \|b\|$,
for some constant $c.$
An {\em $A$-derivation} is a bounded linear map $\delta : A \rightarrow X$
satisfying
$\delta (ab) = a \delta (b) + \delta (a)b$.
An $A$-derivation is {\em inner} if there exists $x$ in $X$ such that
$\delta (a) = ax-xa$.
Finally, $H^1(A,X)$ denotes the quotient of the space of all bounded
derivations by the
space of inner derivations.

\begin{corollary} Let $A$ be a Banach algebra, which has a bounded
faithful representation
as an algebra of operators on a Hilbert space.
Then $H^1(A,X) = 0$ for every bounded $A$-bimodule $X$ if and only if $A$ is
isomorphic
to a finite direct sum of matrix algebras.
\end{corollary}

\begin{remark} Similar results hold for the completely bounded
Hochschild cohomology of an
operator algebra. For the definitions and results see \cite{Pa}.
\end{remark}
\end{document}